\documentclass[10pt]{article}

\usepackage{amsmath,amscd,amsfonts,amsthm}

\newcommand{\shrinkmargins}[1]{
  \addtolength{\textheight}{#1\topmargin}
  \addtolength{\textheight}{#1\topmargin}
  \addtolength{\textwidth}{#1\oddsidemargin}
  \addtolength{\textwidth}{#1\evensidemargin}
  \addtolength{\topmargin}{-#1\topmargin}
  \addtolength{\oddsidemargin}{-#1\oddsidemargin}
  \addtolength{\evensidemargin}{-#1\evensidemargin}
  }

\shrinkmargins{.7}

\DeclareMathOperator{\Hom}{Hom}

\DeclareMathOperator{\PGL}{PGL}

\DeclareMathOperator{\Spec}{Spec}

  \DeclareMathOperator{\Gal}{Gal}

\newcommand{\field}[1]{\mathbb{#1}}

\newcommand{\Q}{\field{Q}}

\newcommand{\Z}{\field{Z}}
\newcommand{\short}{\mathbf{s}}
\newcommand{\F}{\field{F}}

\newcommand{\R}{\field{R}}
\newcommand{\Siegel}{\mathfrak{S}}
\newcommand{\C}{\field{C}}

\newcommand{\A}{\field{A}}
\newcommand{\ord}{\mbox{ord}}

\newcommand{\ra}{\rightarrow}

\newcommand{\OO}{\mathcal{O}}

\newcommand{\Tr}{\mbox{\textnormal Tr}}

\newcommand{\N}{\mbox{\normalfont \bf N}}

\newcommand{\inj}{\hookrightarrow}

\newcommand{\GalK}{\Gal(\bar{K}/K)}

\newcommand{\set}[1]{\{#1\}}
\newcommand{\Disc}{\mathcal{D}}
\newcommand{\DD}{\mathcal{D}}

\newcommand{\beq}{\begin{displaymath}}
\newcommand{\eeq}{\end{displaymath}}
\newcommand{\beqn}{\begin{equation}}
\newcommand{\eeqn}{\end{equation}}
\newcommand{\xx}{\mathbf{x}}
\newcommand{\yy}{\mathbf{y}}

\theoremstyle{plain}
\newtheorem{thm}{Theorem}[section]
\newtheorem{prop}[thm]{Proposition}

\newtheorem{lem}[thm]{Lemma}

\theoremstyle{definition}

\newtheorem{hyp}[thm]{Hypothesis}
\newtheorem{exmp}[thm]{Example}

\theoremstyle{remark}
\newtheorem{rem}[thm]{Remark}

\title{The number of extensions of a number field with fixed degree and bounded discriminant.}
\author{Jordan S. Ellenberg and Akshay Venkatesh\footnote{{\bf  First
    author:}  Department of Mathematics, Princeton University.
  \texttt{ellenber@math.princeton.edu}  Partially supported by NSA
  Young Investigator Grant MDA905-02-1-0097.  {\bf Second author:}
  Department of Mathematics, Massachusetts Institute of Technology.
  \texttt{akshayv@math.mit.edu} Partially supported by NSF Grant DMS-0245606}}
\date{7 September 2003}

\begin{document}
\maketitle

\begin{abstract}
We give an upper bound on the number of extensions of a fixed number field
of prescribed degree 
and discriminant $\leq X$; these bounds improve on work of Schmidt.
We also prove various related results, such as lower bounds for the
number of extensions and upper bounds for Galois extensions.
\end{abstract}

\section{Introduction}

Let $K$ be a number field, and let $N_{K,n}(X)$ be the number of number
fields $L$ (always considered up to $K$-isomorphism) such that $[L:K] = n$ and $\N^K_\Q \DD_{L/K} < X$.  
Here $\DD_{L/K}$ is the relative discriminant of $L/K$,
and $\N^{K}_{\Q}$ is the norm on ideals of $K$, valued in positive integers. 
$\Disc_L = |\DD_{L/\Q}|$ will refer to discriminant over $\Q$.

A folk conjecture, possibly due to Linnik, asserts that
\beq
N_{K,n}(X) \sim c_{K,n} X \ \ \ \ \ \ \ \   (n \mbox{ fixed, }X \rightarrow \infty)
\eeq

This conjecture is trivial when $n=2$; it has been proved for $n=3$ by
Davenport and Heilbronn~\cite{dave:dahe} in case $K=\Q$, and by
Datskovsky and Wright in general~\cite{dats:cubics}; and for $n=4,5$
and $K=\Q$ by Bhargava~\cite{bhar:thesis},\cite{bhar:hcl4}.  
A weaker
version of the conjecture for $n=4$ and arbitrary $K$ 
is due to Yukie.
These beautiful results are proved
by methods which seem not to extend to higher $n$.  The best upper
bound for general $n$ is due to Schmidt~\cite{schm:columbia}, who
showed
\beq
N_{K,n}(X) \ll X^{(n+2)/4}
\eeq
where the implied constant depends on $K$ and $n$.
We refer to \cite{cohen} for a survey of results. 

In many cases, it is easy to show that $N_{K,n}(X)$ is bounded below
by a constant multiple of $X$; for instance, if $n$ is even, simply
consider the set of quadratic extensions of a fixed $L_0/K$ of degree
$n/2$.  For the study of lower bounds it is therefore more interesting
to study the number of number fields $L$ such that 
$[L:K] = n$, $\N^K_\Q \DD_{L/K}
< X$ and the Galois closure of $L$ has Galois group $S_n$ over $K$.  Denote this number by $N'_{K,n}(X)$.
Malle showed~\cite[Prop.\ 6.2]{mall:jnt} that 
\beq
N'_{\Q,n}(X) > c'_n X^{1/n}
\eeq
for some constant $c'_n$.

The main result of this paper is to improve these bounds, with
particular attention to the ``large $n$ limit.'' The upper
bound lies much deeper than the lower bound. 

Throughout this paper we will use $\ll$ and $\gg$ 
where the implicit constant depends on $n$; we will not make this
$n$-dependency explicit. 

\begin{thm}
For all $n>2$ and all number fields $K$, we have 
\beq
N_{K,n}(X) \ll (X \DD_K^n A_n^{[K:\Q]})^{\exp(C
  \sqrt{\log n})}
\eeq
where $A_n$ is a constant depending only on $n$, and $C$ is an absolute
constant. Further,
\beq
X^{1/2 + 1/n^2} \ll_K N'_{K,n}(X)
\eeq
% that can be taken to be $(3/2)(1 + \log 4)$.
\label{th:main}
\end{thm}

In particular, for all $\epsilon > 0$
\begin{equation} \label{eqn:loglimits}
\limsup_{X \ra \infty} \frac{\log N_{K,n}(X)}{\log X}
\ll_{\epsilon} n^{\epsilon}, \ \ \liminf_{X \ra \infty} \frac{\log
N_{K,n}'(X)}{\log X} \geq \frac{1}{2} + \frac{1}{n^2}
%48 \exp (C \sqrt {\log n}).
\end{equation}
Linnik's conjecture claims that the limit in (\ref{eqn:loglimits}) is
equal to $1$; thus, despite its evident imprecision, the upper bound
in Theorem~\ref{th:main} seems to offer the first serious evidence
towards this conjecture for large $n$. It is also worth observing that
de Jong and Katz~\cite{katz:mgpoints} have studied a problem of a
related nature where the number field $K$ is replaced by the function
field $\F_q(T)$; even here, where much stronger geometric techniques
are available, they obtain an exponent of the nature $c \log(n)$.
This suggests that replacing $n^\epsilon$ in \eqref{eqn:loglimits} by
a constant will be rather difficult.

\medskip

We will also prove various related results on the number of number
fields with certain Galois-theoretic properties.  For instance, if $G
\leq S_n$, let $N_{K,n}(X;G)$ be the number of number fields $L$ such
that $[L:K] = n$, $\mathbf{N}^K_{\Q} \mathcal{D}_{L/K} < X$, and the
action of $\GalK$ on embeddings $K \inj \C$ is conjugate to the
$G$-action on $\{1, \dots, n\}$.  We describe how one can
obtain upper bounds on $N_{K,n}(X;G)$ using the invariant theory of
$G$.  A typical example is:

\begin{prop}
Let $G \leq S_6$ be a permutation group whose action is conjugate to the
$\mathrm{PSL}_2(\mathbb{F}_5)$-action
on $\mathbb{P}^1(\mathbb{F}_5)$. 
Then $N_{\Q,6}(X;G) \ll_{\epsilon} X^{8/5+\epsilon}$.
\end{prop}

Specializing further, let $N_{K,n}(X; Gal)$
be the number of {\em Galois} extensions among those counted by
$N_{K,n}(X)$; we prove the following upper bound.

\begin{prop}
For each $n > 2$, one has $N_{K,n}(X; Gal) \ll_{K,n,\epsilon}
X^{3/8+\epsilon}$.  
\end{prop}

In combination with the lower bound in Theorem~\ref{th:main}, this
shows that if one orders the number fields of fixed degree over $\Q$
by discriminant, a random one is not Galois. 

Although we will use certain {\em ad hoc} tools, the central idea will
always be to count fields by counting integral points on certain
associated varieties, which are related to the invariant theory of the
Galois group. These varieties must be well-chosen to obtain good
bounds. In fact, the varieties we use are birational to 
the Hilbert scheme of $r$ points in $\mathbb{P}^n$, suggesting
the importance of a closer study of the distribution of rational points
on these Hilbert schemes. 

The results can perhaps be improved using certain techniques
from the study of integral points, such as the result of Bombieri-Pila
\cite{pila}.  However, the proof of Theorem 1.1 turns out, somewhat
surprisingly, to require only elementary arguments from the geometry
of numbers and linear algebra.

\section{Proof of Upper Bound} 
The main idea of Schmidt's proof is as follows: by Minkowski's
theorem, an extension $L/K$ contains an integer $\alpha$ whose archimedean
valuations are all bounded by a function of $\Delta_L = \N^K_\Q \DD_{L/K}$.
Since all the archimedean absolute values are bounded in terms of
$\Delta_L$, so are the symmetric functions of
these absolute values;  in other words, $\alpha$ is a root of a monic
polynomial in $\Z[x]$ whose coefficients have (real) absolute value
bounded in terms of $\Delta_L$.  There are only finitely many such
polynomials, and counting them gives the theorem of
\cite{schm:columbia}.

The main idea of Theorem~\ref{th:main} is to count $r$-tuples of
integers in $L$ instead of single integers. 

Let $\A^n = \Spec(\Z[x_1,x_2, \dots, x_n])$ denote affine $n$-space,
which we regard as being defined over $\Z$.  We fix an algebraic
closure $\bar{K}$ of $K$.  Let $\rho_1, \ldots, \rho_n$ be the
embeddings of $L$ into $\bar{K}$.  Then the map $\phi_L = \rho_1
\oplus \ldots \oplus \rho_n$ embeds $\OO_L$ in $\bar{K}^n =
\A^n(\bar{K}) $, and the direct sum of $r$ copies of this map is an
embedding $\OO_L^r \ra (\bar{K}^n)^r = (\A^n)^r(\bar{K})$ (which we
also, by abuse of notation, call $\phi_L$.)

The affine variety $(\A^n)^r$ is naturally coordinatized by functions
$\set{x_{j,k}}_{1 \leq j \leq n, 1 \leq k \leq r}.$ The symmetric
group $S_n$ acts on $(\A^n)^r$ by permuting $x_{1,k}, \ldots,
x_{n,k}$ for each $k$. The $S_n$-invariants in the coordinate ring of
$(\A^n)^r$ are called {\em multisymmetric functions}.  If $f$ is a
multisymmetric function, the composition $f \circ \phi_L: \OO_L^r \ra
\bar{K}$ takes image in $\OO_K$.  It follows that if 
$R \subset \Z[\set{x_{j,k}}_{1 \leq j \leq n, 1 \leq k \leq r}]^{S_n}$
is a subring of the ring of multisymmetric functions and $A = \Spec(R)$, 
there is a map of sets
\beq
F:  \bigcup_{L} \OO_L^r \ra A(\OO_K)
\eeq
where the union is over all number fields $L$ with $[L:K] = n$. 

Our overall strategy can now be outlined as follows.  If $x$ is an algebraic integer, write $||x||$ for the maximum of the archimedean absolute values of $x$.  For a positive real
number $Y$, let $B(Y)$ be the set of algebraic integers $x$
in $\overline{K}$ with degree $n$ over $K$ and $||x|| < Y$. Let $f_1, \dots, f_s
\in \Z[\set{x_{j,k}}_{1 \leq j \leq n, 1 \leq k \leq r}]^{S_n}$
be multisymmetric functions
with degrees $d_1, \dots, d_s$.
Put $R = \Z[f_1, \dots, f_s]$, and set $A = \Spec(R)$. 
Then there is a constant
$c$ such that 
(for any $Y$) one has $||f_i(\phi_L(\alpha_1,\alpha_2, \dots, \alpha_r))|| < c Y^{d_i}$ 
whenever $\alpha_j \in B(Y) \, (1 \leq j \leq r)$.   Let $A(\OO_K)_Y$ be the subset of $A(\OO_K)$
consisting of points $P$ such that $||f_i(P)|| < c Y^{d_i}$.  Then for any
subset $S_Y$ of $B(Y)^r$, we have 
a diagram of sets
\begin{equation}
\begin{CD}
\set{(L,\mathbf{\alpha}_1,
\mathbf{\alpha}_2, \dots, \mathbf{\alpha}_r):[L:K] = n, \Delta_L < X, (\mathbf{\alpha}_1, \ldots, \mathbf{\alpha}_r) \in
  (\OO_L)^r \cap S_Y} @>F>> A(\OO_K)_Y \\
@VVV \\
\set{L: [L:K], \Delta_L < X}
\end{CD}
\label{eq:setdiag}
\end{equation}
The cardinality of the lower set is precisely $N_{K,n}(X)$.Our goal is to choose $A = \Spec(R), Y,$ and $S_Y$ in such a way that that the
vertical map in \eqref{eq:setdiag} is surjective (by Minkowski's
theorem), while the horizontal map $F$ has finite fibers whose
cardinality we can bound.  This will yield the desired bound on
$N_{K,n}(X)$.  Since $|A(\OO_K)_Y| \ll_K (c^s
Y^{\sum_i d_i})^{[K:\Q]}$, it should be our aim to choose $f_1, \dots, f_s$
whose total degree is as low as possible.

%Our strategy will be to construct such a variety $A$ which is large enough
%so that $\phi$ is ``close to injective'', but which is small enough so
%that we can put a strong bound on the number of points in $R(\OO_K)$
%with bounded height.  This will yield the desired bound on the number
%of fields $L$ which contribute to the left hand side.

We begin with a series of lemmas about polynomials over an arbitrary
characteristic-$0$ field $F$.

%(we will
%need statements that guarantee that $\phi$ is close to injective, but they
%can be checked ``geometrically,'' i.e. over $k$.)

Let $S$ be any test ring. 
We give $\A^n$ the structure of a ring scheme
so that the ring structure on $\A^n(S) = S^n$ is the natural one.
Let $\Tr$ be the map $\A^n \rightarrow \A^1$
which, on $S$-points, induces the map $(z_1,
\ldots, z_n) \in S^n\mapsto z_1 + \ldots + z_n \in S$.
Given an element $\xx =(x_{j,k})_{1 \leq j \leq n, 1 \leq k \leq r} 
\in (S^n)^r$, we denote by $\xx_k \in S^n$ the $k$th ``row''
$(x_{1,k},x_{2,k},\dots, x_{n,k})$, and by $\xx^{(j)} \in S^r$
the $j$th ``column'' $(x_{j,1}, x_{j,2}, \dots, x_{j,r})$. 
These correspond to maps $\xx \mapsto \xx_k: (\A^n)^r \rightarrow \A^n,
\xx \mapsto \xx^{(j)}: (\A^n)^r \rightarrow \A^r$. 

Let
$\sigma = (i_1, \ldots, i_r)$ be an element of $\Z_{\geq 0}^r$;
we will think of $\Z_{\geq 0}^r$ as an additive semigroup, operations
being defined pointwise.
Then
$\sigma$ defines a $S_n$-equivariant
map $\chi_\sigma: (\A^n)^r \ra \A^n$ by the rule
\beq
\chi_\sigma(\xx) = \xx_1^{i_1} \xx_2^{i_2} \dots \xx_r^{i_r}
\eeq

In particular, $F^n = \A^n(F)$ has a ring structure,
and $\Tr, \chi_{\sigma}$ induce maps on $F$-points, namely $\Tr: F^n  \rightarrow F,
\chi_{\sigma}: (F^n)^r \rightarrow F^n$; we abuse notation and use the same
symbols for these maps. 
The map $(x,y) \mapsto \Tr(xy)$
 is a nondegenerate pairing on $F^n$, with respect to which we can speak of
 ``orthogonal complement''.

\begin{lem} Let $\xx \in (F^n)^r$, and let $\Sigma$ be a
  subset of $\Z_{\geq 0}^r$ such that the $|\Sigma|$ vectors
  ${\chi_\sigma(\xx)}_{\sigma \in \Sigma}$ generate a subspace of $F^n$
(considered as an $F$-vector space)
  of dimension greater than $n/2$.  

  Denote by $\Sigma + \Sigma$ the set of sums of two elements of
  $\Sigma$.  Let $W \subset F^n$ be the subspace of $F^n$ spanned by
  ${\chi_\sigma(\xx)}_{\sigma \in (\Sigma + \Sigma)}$.  Then the
  orthogonal complement of $W$ is contained in a coordinate hyperplane
  $x_j = 0$ for some $j$. 
\label{le:n2}
\end{lem}

\begin{proof}
Write $m$ for $|\Sigma|$ and let $v_1, \ldots, v_m$ be the vectors
$\chi_\sigma(\xx)$ as $\sigma$ ranges over $\Sigma$.  Then $W$ is the
space spanned by the products $v_a v_b$ (the algebra structure on $F^n$
structure being as noted above). Suppose $w$
is orthogonal to $W$. Then 
\begin{equation}
\Tr (v_a v_b w) = 0
\label{eq:tr0}
\end{equation}
for all $a,b$; if $V$ is the space spanned by the $\set{v_a}$, then
\eqref{eq:tr0} implies that $wV$ and $V$ are orthogonal.  This implies
in turn that $\dim wV \leq n - \dim V < \dim V$, so multiplication by
$w$ is not an automorphism of $F^n$; in other words, $w$ lies on a
coordinate hyperplane.  A subspace of $F^n$ contained in a union of
coordinate hyperplanes is contained in a single coordinate hyperplane;
this completes the proof.
\end{proof}

For each $\sigma \in
\Z_{\geq 0}^r$, let $f_\sigma: (\A^n)^r \rightarrow \A^1$ be the composition $\Tr \circ
\chi_\sigma$.  Then $f_\sigma$ is a multisymmetric function.
When $\Sigma$
is a subset of $\Z_{\geq 0}^r$, we denote by $R_\Sigma$ the subring
of functions on $(\A^n)^r$
generated by $\set{f_\sigma}_{\sigma \in \Sigma}$.  
One has a natural map of affine schemes
\begin{equation}
F_\Sigma:  (\A^n)^r \ra \Spec R_\Sigma.
\label{eqn:fsigma} \end{equation}
The goal of the algebro-geometric part of our argument is to show that,
by choosing $\Sigma$ large enough, we can guarantee that $F_\Sigma$ is
generically finite, and even place some restrictions on the locus in
$(\A^n)^r$ where $F_\Sigma$ has positive-dimensional fibers.

\begin{lem}  Let $\xx$ be a point of $(\A^n)^r(F)$, and let $\Sigma$ be a
  subset of $\Z_{\geq 0}^r$ such that the $|\Sigma|$ vectors
  ${\chi_\sigma(\xx)}_{\sigma \in \Sigma}$ span $F^n$ as an $F$-vector
  space.  For each $k$
  between $1$ and $r$ let $e_k \in \Z_{\geq 0}^r$ be the vector with a
  $1$ in the $k$'th coordinate and $0$'s elsewhere.
  Let $\Sigma'$ be a set which contains $\Sigma + \Sigma$, and $\Sigma
  + e_k$ for all $k$.

%  Then $F_{\Sigma'}$ is finite on an open neighborhood of
%  $\xx$.   
Then the preimage $F_{\Sigma'}^{-1}(F_{\Sigma'}(\xx))
\subset (\A^n)^r(F)$ is finite, of
cardinality at most $(n!)^r$. 
\label{le:fsigfin}
\end{lem}

\begin{proof}  Let $\yy$ be $F_{\Sigma'}(\xx)$.  Let $m = |\Sigma|$. 
%It suffices to show
%  that there are only finitely many $\xx'$ such that $F_{\Sigma'}(\xx') =
%  \yy$.  

  As in the proof of Lemma~\ref{le:n2}, let $v'_1, \ldots, v'_m$ be
  the image of $\xx$ under the $\set{\chi_\sigma}_{\sigma \in
  \Sigma}$. We may suppose by relabeling that
  $v'_1, \ldots, v'_n$ form a basis for $F^n$ (as an $F$-vector space).  Since $\Sigma'$
  contains $\Sigma + \Sigma$, the determination of $\yy$ fixes
  $\Tr(v'_a v'_b)$ for all $a,b$; and since $\Sigma'$ contains $\Sigma
  + e_k$, we also know the traces $\Tr(v'_a \xx_{k})$ for
  all $a$ and $k$. It follows that, for each $k$, we can represent the
  action of multiplication by $\xx_{k}$ on the $F$-vector space spanned by
  $v'_1, \ldots, v'_n$ by a matrix whose coefficients are determined
  by $\yy$.  But such a matrix evidently determines $\xx_{k}$ up to permutation of coordinates; this proves the desired
  result.
\end{proof}

% The above lemmas allow us to show, for particular choices of $\Sigma$,
% that the ``bad locus'' of $F_\Sigma$ does not contain certain linear
% subspaces.  

In the proof of Proposition~\ref{prop:rtuple} below, we will need to show
that, by allowing $\xx$ to vary over certain subspaces of $(F^n)^r$,
we can ensure that $\xx$ can be chosen in order to verify the
hypothesis of Lemma~\ref{le:n2}.
   
\begin{lem} Let $V$ be a $F$-subspace of $F^n$ of dimension $m$,
  and let $\Sigma \subset \Z_{\geq 0}^r$ be a subset of size $m$.
  Let $Z \subset V^r$ be the subset of points $\xx \in V^r$
  such that  the $m$ vectors $\chi_\sigma(\xx)_{\sigma \in \Sigma}$
  are not linearly independent (over $F$) in $F^n$.  Then $Z$ is not the whole of
  $V^r$. If one identifies $V^r$ with $F^{mr}$,
  $Z$ is contained in the $F$-points of 
  a hypersurface, defined over $F$, whose degree 
  is bounded by a
  constant depending only on $n$ and $\Sigma$.
\label{le:indep}
\end{lem}

\begin{proof}
We may assume (by permuting coordinates) that the
map ``projection onto the first $m$ coordinates,'' which we denote
$\pi:F^n \ra F^m$, induces an isomorphism $V \cong F^m$.  Suppose there is a nontrivial linear relation
\begin{equation} 
\sum_{\sigma \in \Sigma} c_\sigma \chi_\sigma(\xx) = \mathbf{0} \in F^n,
\label{eq:relation}
\end{equation}
that is, suppose $\xx \in Z$.  Each $\sigma \in \Sigma$ also defines a
map $F^r \rightarrow F$ (derived from the map $\chi_{\sigma}:
(\A^n)^r \rightarrow \A^n$ with $n=1$) so we may speak of $\chi_{\sigma}(\xx^{(j)})
\in F$ for $1 \leq j \leq n$.  By abuse of notation we also use $\pi$
to denote the projection of $(F^{n})^r$ onto $(F^{m})^r$.  Then the
restriction of $\pi$ to $V^r$ is an isomorphism $V^r \cong F^{mr}$.

%The maV^r \rightarrow F^{mr}$ given by $\xx \mapsto (\xx^{(1)}, \xx^{(2)},
%\dots, \xx^{(m)})$ is an isomorphism,

Any nontrivial linear relation between the $\chi_\sigma(\xx)$ yields a
nontrivial relation between the $m$ vectors $\chi_\sigma(\pi(\xx))$ in
$F^m$. This in turn implies vanishing of the determinant
\beq
D = \left| \begin{array}{ccc}
\chi_{\sigma_1} (\xx^{(1)}) & \cdots & \chi_{\sigma_1}(\xx^{(m)}) \\
\vdots & & \vdots \\
\chi_{\sigma_m} (\xx^{(1)}) & \cdots & \chi_{\sigma_m}(\xx^{(m)})
\end{array} \right|.
\eeq

The contribution of each $m \times m$ permutation matrix to $D$ is a
distinct monomial in the $mr$ variables, so $D$ is not identically $0$
in $F[x_{1,1}, \ldots, x_{m,r}]$.  Evidently the degree of $D$ is
bounded in terms of $n$ and $\Sigma$.  Let $V(D)$ be the vanishing
locus of $D$ in $(F^m)^r$. 
Now the locus in $Z$ is
contained in $\pi^{-1}(V(D))$, which yields the desired result.
\end{proof}

%  I find this the most confusing of the lemmas; Akshay, can you see a
%  way to explain it more clearly?

Finally, we need a straightforward fact about points of low height on
the complements of hypersurfaces.

\begin{lem}
Let $f$ be a polynomial of degree $d$ in variables $x_1, \ldots,
x_n$.  Then there exist integers $a_1, \ldots, a_n$ such that $\max_{1 \leq
i \leq n}
|a_i| \leq (1/2)(d+1)$ and $f(a_1, \ldots, a_n) \neq 0$.
\label{le:hypercomp}
\end{lem}

\begin{proof}  There are at most $d$ hyperplanes
  on which $f$ vanishes, which means that the function $g(x_2, \ldots,
  x_n) := f(a_1, x_2, \ldots, x_n)$ is not identically $0$ for some
  $a_1$ with absolute value at most $(1/2)(d+1)$.  Now proceed by
  induction on $n$.
\end{proof}

Now we are ready for the key point in the proof of
Theorem~\ref{th:main}.  The point is to use the lemmas above to
construct $\Sigma$ which is small enough that $\Spec R_\Sigma$ has few
rational points of small height, but which is large enough so that
$F_\Sigma$ does not have too many positive-dimensional fibers.

\begin{prop} \label{prop:rtuple}Let $\Sigma_0$ be a subset of $\Z_{\geq 0}^r$ of size $m
  > n/2$; let $\Sigma_1 \subset \Z^r_{\geq 0}$ contain $\Sigma_0 + \Sigma_0$;
    and let $\Sigma \subset \Z^r_{\geq 0}$ contain $\Sigma_1 + \Sigma_1$ and
    $\Sigma_1 + e_k$ for all $k$.

Let $L$ be a finite extension of $K$ with $[L:K] = n$.  Then there is
an $r$-tuple $(\alpha_1, \ldots, \alpha_r) \in \OO_L^r$ such that
\begin{itemize}
\item For every $k$, we
  have
\beq
||\alpha_k|| \ll_\Sigma \DD_L^{1/d(n-2)},
\eeq
where $d = [K:\Q]$.

\item The set $F_\Sigma^{-1} (F_\Sigma( (\phi_L(\alpha_1, \ldots, \alpha_r))))
  \subset (\A^n)^r(\bar{K})$ has cardinality at most $(n!)^r$.

\item The elements $\alpha_1, \ldots, \alpha_r$ generate the field
  extension $L/K$.

\end{itemize}
\end{prop}

\begin{proof}

First of all, note that if $(\alpha_1, \ldots, \alpha_r),\Sigma_0,
\Sigma_1, \Sigma$ satisfy the conditions above, then so do $(\alpha_1,
\ldots, \alpha_r),\Sigma'_0, \Sigma_1, \Sigma$ for any subset
$\Sigma'_0 \subset \Sigma_0$.  So it suffices to prove the theorem in
case $n/2 < m \leq (n/2 + 1)$.

%For every $\beta \in \OO_L$, define $||\beta||$ to be the maximal archimedean absolute value of $\beta$.
%
%  Revised up to here 16 Jun

Let $1=\beta_1, \ldots, \beta_{nd}$ be a $\Q$-linearly independent set
of integers in $\OO_L$ such that $||\beta_i||$ is the $i$'th successive
minimum of $||\cdot||$ on $\OO_L$, in the sense of Minkowski's second
theorem~\cite[III, \S 3]{sieg:gon}.  The $K$-vector space spanned by
$\beta_1, \ldots, \beta_{md}$ has $K$-dimension at least $m$, so we may
choose $\gamma_1, \ldots, \gamma_m$ among the $\beta_i$ which are
linearly independent over $K$.

Let $V \subset \bar{K}^n$ be the $\bar{K}$-vector space spanned by $\set{\phi_L(\gamma_i)}_{1
\leq i \leq m}$. Then by Lemma~\ref{le:indep} there is a constant
$C_{n,\Sigma_0}$ and a hypersurface $Z \subset V^r$ of degree
$C_{n,\Sigma_0}$ such that, for all $\xx$ not in $Z(\bar{K})$, the $m$ vectors
${\chi_\sigma(\xx)}_{\sigma \in \Sigma_0}$ are $\bar{K}$-linearly
independent in $\bar{K}^n$.

For every field $M$ strictly intermediate between $K$ and $L$, we let $V_M \in
\bar{K}^n$ be the $\bar{K}$-vector subspace $\phi_L(M) \subset \bar{K}^n$.  
Each $(V \cap V_M)^r$ is a certain linear subspace of $V^r$;
note that, since $m > n/2$, no subspace
$V_M$ contains $V$. Let
$Z'$ be the union of $Z(\bar{K})$ with $(V \cap V_M)^r$, as $M$ ranges over all
fields between $K$ and $L$.  

%  Up to here, Jul 22

Now let $Y$ be a hypersurface of $V^{r}$ so that $Y(\bar{K})$ contains $Z'$; 
one may choose $Y$ so that the degree of $Y$
is bounded in terms of $n$ and $\Sigma_0$.  By
Lemma~\ref{le:hypercomp}, there is a constant $H$, depending only on
$n$ and $\Sigma_0$, so that, for any lattice $\iota: \Z^{mr} \inj V^{r}$
(i.e. we require $\iota(\Z^{mr})$ spans $V^r$ over $\bar{K}$)
there is a point $p \in \Z^{mr}$, with $\iota(p) \notin Y(\bar{K})$,  whose
coordinates have absolute value at most $H$.

It follows that there exists a set of $mr$ integers $c_{1,1}, \ldots,
c_{m,r}$ with $|c_{j,k}| \leq H$, such that
\beq
\xx = (\phi_L(c_{1,1} \gamma_1 + \ldots + c_{m,1} \gamma_m),
\ldots, \phi_L(c_{1,r} \gamma_1 + \ldots + c_{m,r} \gamma_m))
\eeq
is not in $Y(\bar{K})$.  For each $k$ between $1$ and $r$ define $\alpha_k \in
\mathcal{O}_L$ via
\beq
\alpha_k = c_{1,k} \gamma_1 + \ldots + c_{m,k} \gamma_m.
\eeq

Let $W \subset L$ be the $K$-subspace spanned by
$\chi_\sigma(\alpha_1, \ldots, \alpha_r)$ as $\sigma$ ranges over
$\Sigma_1$ (here we regard $\chi_{\sigma}$ as a map $L^r \rightarrow L$,
c.f. remarks after (\ref{eq:relation})). Suppose $W$ is not the whole of $L$.  Then there is a
nonzero element $t \in L$ such that $\Tr^L_K tw = 0$ for all $t \in
W$.  It follows that $\phi_L(t) \in \bar{K}$ lies in the orthogonal complement (w.r.t to the form $\Tr$ on $\bar{K}^n$) of
$\phi_L(W) \subset \bar{K}^n$.  But the orthogonal complement to the $\bar{K}$-span of
$\phi_L(W)$ is contained in a coordinate hyperplane by
Lemma~\ref{le:n2}.  Since $\rho_j(t)$ cannot be $0$ for any $j$ and
any nonzero $t$, this is a contradiction; we conclude that $W = L$,
and thus that the vectors $\set{\chi_\sigma(\xx)}_{\sigma \in
\Sigma_1}$ span $L$ as a $K$-vector space. 

The bound
on the size of the fiber $F_\Sigma^{-1}(F_\Sigma(\xx))$ follows
from Lemma~\ref{le:fsigfin},
and
the fact that $\xx \notin V_M^r$ for any $M$ implies that $\alpha_1,
\ldots, \alpha_r$ generate the extension $L/K$.
 
It remains to bound the archimedean absolute values of the
$\alpha_i$. 
The image of $\OO_L$ in $\OO_L \otimes_{\Z} \R$ is
a lattice of covolume $\DD_L^{1/2}$, so by Minkowski's
second theorem\cite[Thm 16]{sieg:gon}, 
\beq
\prod_{i=1}^{nd} ||\beta_i|| \leq  \DD_L^{1/2}.
\eeq
The $||\beta_i||$ form
a nondecreasing sequence, so for $m < n$, we have
\beq
||\beta_{md}||^{(n-m)d} \leq \prod_{i=md+1}^{nd} ||\beta_i|| \leq
 \DD_L^{1/2}.  
\eeq
Since $m \leq (1/2)n + 1$, we get
\beq
||\beta_i|| <  (\DD_L)^{1/d(n-2)}
\eeq
for all $i \leq m$. It follows that all archimedean absolute values
of $\gamma_i$ for $i \leq m$ are bounded by a constant multiple of
$\DD_L^{1/d(n-2)}$, the implicit constant being absolute.   The result follows, since each $\alpha_k$ is an
integral linear combination of the $\gamma_i$ with coefficients bounded by
$H$.
\end{proof}

We are now ready to prove the upper bound in Theorem~\ref{th:main};
what remains is merely to make a good choice of $\Sigma$ and apply
Proposition~\ref{prop:rtuple}. Let $r$ and $c$ be positive integers
such that ${r + c \choose r} > n/2$, and
let $\Sigma_0$ be the set of all $r$-tuples of nonnegative integers with sum
at most $c$.  We shall choose $r,c$ in the end; but
$r,c,\Sigma_0,\Sigma$ will all depend only on $n$,
so that all constants that depend on them in fact depend only on $n$. 

%%Then $|\Sigma_0| = {r+c \choose r}$.  We choose $c$ to be
%jithe smallest integer such that ${r+c \choose r} > n/2$.  We assume
%from the start that $r$ and $c$ are less than $n$, so that all
%constants depending on $\Sigma$ depend only on $n$.

%Let $\delta > 0$ be a constant, let $r$ be
%the greatest integer less than $(1+\delta)\sqrt{\log n}$, and let
%$c$ be the greatest integer less than $e^\sqrt{\log n}$.  Let
%$\Sigma$ be the set of all $r$-tuples of nonnegative integers with sum
%at most $c$.  Then $|\Sigma_0| = {r+c \choose r}$; by Stirling's
%approximation,
%\beq
%\log |\Sigma_0| \sim r \log c = (1+\delta) \log n
%\eeq
%as $n \ra \infty$, so $|\Sigma_0| > n/2$ for $n$ large enough.  

Now take $\Sigma$ to be the set of all $r$-tuples of nonnegative
integers with sum at most $4c$, and consider the map
\beq
F_\Sigma:  (\A^n)^r \ra \Spec R_\Sigma.
\eeq
By Proposition~\ref{prop:rtuple}, to every field $L$ with $[L:K]=n$ we can
associate an $r$-tuple $(\alpha_1, \ldots, \alpha_r)$ of integers
satisfying the three conditions in the statement of the proposition.  Define
$Q_L \in (\A^n)^r(\bar{K})$ to be $\phi_L(\alpha_1, \ldots, \alpha_r)$, and
let $P_L \in \Spec R_\Sigma(\OO_K)$ be the point $F_\Sigma(Q_L)$.

By the second condition on $\alpha_1, \ldots, \alpha_r$, there are at
most $(n!)^r$ points in $F_\Sigma^{-1}(P_L)$.  By the third condition,
$Q_L = Q_{L'}$ only if $L$ and $L'$ are isomorphic over $K$. 
We conclude that at most
$(n!)^r$ fields $L$ are sent to the same point in $\Spec
R_\Sigma(\OO_K)$.

We now restrict our attention to those fields $L$ satisfying
\beq
\N^K_\Q \DD_{L/K} < X.
\eeq
In this case, for every archimedean valuation $|\cdot|$ of $L$ and
every $k \leq r$ we have the bound
\begin{equation}
|\alpha_k| \ll \DD_L^{1/d(n-2)} \ll (X \DD_K^n)^{1/d(n-2)}
\label{eq:cav}
\end{equation}
 
Now, $f_{\sigma}$ being as defined prior to (\ref{eqn:fsigma}),  $f_\sigma(Q_L)$ is an element of $\OO_K$, which (by choice of $\Sigma$)
we can express as a polynomial of degree at most $4c$ (and absolutely
bounded coefficients) in the numbers $\rho_j(\alpha_k) \in \bar{K}$.
If $|\cdot|$ is any archimedean absolute value on $K$, we can extend $|\cdot|$ to a
archimedean absolute value on $L$, and by \eqref{eq:cav} we have
\beq
|f_\sigma(Q_L)| \ll (X \DD_K^n)^{4c/d(n-2)}.
\eeq

%Since this argument depends in no way the original choice of embedding
%$\rho: K \ra \C$, we may conclude that $f_\sigma(Q_L)$ (considered as
%an algebraic integer in $K$) has {\em all} archimedean absolute values
%bounded by $B = A_2 (X \DD_{K/\Q}^n)^{4c/d(n-2)}$.  

The number of elements of $\OO_K$ with archimedean absolute values at
most $B$ is $ \leq (2B+1)^d$.  (For large enough $B$, one can save an extra
factor of $\DD_K^{1/2}$; this is not necessary for our purpose.)
In view of the above equation, the number of possibilities for
$f_{\sigma}(Q_L)$ is  $\ll (X \DD_K^n A_n^d)^{4c / (n-2)}$
where $A_n$ is a constant depending only on $n$.

% We conclude that the number of
%possibilities for $f_\sigma(Q_L)$ is at most
%\beq
%A_4 X^{4c/(n-2)} \D_{K/\Q}^{4cn/(n-2)}.
%\eeq
%(wE DISCard the insignificant $\DD_{K/\Q}^{-1/2}$.)

Now the point $P_L \in \Spec R_\Sigma(\OO_K)$ is determined by 
 $f_\sigma(Q_L)$ ($\sigma \in \Sigma$) and we have 
$
|\Sigma| = {r + 4c \choose r}$.  
The number of possibilities for $P_L$ is therefore 
%So the number of points $P$ of $\Spec R_\Sigma(\OO_K)$ such
%that $f_\sigma(P)$, for all $\sigma \in \Sigma$, has all archimedean
%absolute values bounded by $B$ is at most
$\ll 
(X \DD_K^n A_n^d)^{(4c/(n-2))  {r + 4 c \choose r}}$. 

Since each number field $L$ contributes a point to this count, and
since no point is counted more than $(n!)^r$ times, we have
\begin{equation}
N_{K,n}(X) \ll   (X \DD_K^n A_n^d)^{(4c/(n-2))  {r + 4 c
    \choose r}}. \label{eq:nknrc}
\end{equation}

Now is a suitable time to optimize $r$ and $c$. 
We may assume $n \geq 3$. 
Take $r$ to be the greatest integer $\leq \sqrt{\log(n)}$,
and choose $c$ to be the least integer $\geq (n r!)^{1/r}$.
Note that $c \geq n^{1/r} \geq e^{\sqrt{\log(n)}} \geq e^r  \geq r$
and $c \leq 2 (nr!)^{1/r}$. 
Then ${r + c \choose r} > c^r/r! > n$
whereas ${r + 4c \choose r} \leq {5 c \choose r}
\leq \frac{(5c)^r}{r!} \leq 10^r n$. 
%$ minimal such that ${r + c \choose r} > n/2$;
%then f we choose $r$ on
%the order of $\sqrt(\log n)$, then $c$ is much larger than $r$, and
%${r + c \choose r}$ is on order of $c^r$.  We conclude that $c$ is on
%the order of $n^{1/r}$, which is to say $c$ is on the order of
%$\exp(\sqrt(\log n))$. 
Substituting these values of $r,c$ into
\eqref{eq:nknrc} yields the upper bound of Theorem~\ref{th:main}.

%  Take $r$ to be the least
%integer greater than $\sqrt{\log n}$.  Since 
%\beq
%\frac{c}{r+c} {r+c \choose r} = {r+c-1 \choose r} 
%\leq n/2 
%\eeq
%we have 
%\beq
%{r + c \choose r} \leq (1/2)(1+r/c) n
%\eeq
%It is easy to see that $c > r$ for $n$ large enough, so 
%\beq
%{r + c \choose r} \leq n.
%\eeq
%In order to get an upper bound for $c$, write
%\beq
%n \geq {r + c \choose r} \geq (c/r)^r
%\eeq
%whence $c \leq rn^{1/r}$.  Since the function $rn^{1/r}$ is decreasing
%in $r$ when $r > \log n$, we have
%\beq
%c \leq (\sqrt{\log n}( n^{1/\sqrt{\log n}} = (\sqrt{\log n}) e^{\sqrt{\log
%    n}}.
%\eeq

%Plugging the upper bounds for $r,c,$ and ${r+c \choose c}$ into
%\eqref{eq:nknrc} now yields
%\beq
%N_{K,n}(X) \ll   (X \DD_{K/\Q}^n)^{(4n\sqrt{\log n} \exp(\sqrt{\log
%    n}/(n-2)) 4^{\sqrt{\log n} + 1}}
%\eeq
%or, more simply, for $n > 2$, 
%\beq
%N_{K,n}(X) \ll  (X \DD_{K/\Q}^n)^{\exp(C
%  \sqrt{\log n})}
%\eeq
%where $C$ is an absolute constant.  This is the upper bound of Theorem~\ref{th:main}.

In the language of the beginning of this section, we have taken $A$ to
be $\Spec R_\Sigma$, the map $F$ to be $F_\Sigma$, and the set $S_Y$
to be the set of $r$-tuples of integers $\alpha_1, \ldots, \alpha_r$ of
so that $\alpha_j \in B(Y) (1 \leq j \leq r)$
and whose image under $\phi_L$ lies in $V^r - Z'$
(notation of proof of Proposition~\ref{prop:rtuple}.)  Minkowski's theorem guarantees that
each number field $L$ contains an $r$-tuple of integers in $S_Y$ for
some reasonably small $Y$, while the lemmas leading up to
Proposition~\ref{prop:rtuple} show that the fibers of $F$ containing a
point of $S_Y$ have cardinality at most $(n!)^r$.

Another way to think of the method is as follows:  we can factor
$F_\Sigma$ as 
\beq
(\A^n)^r \ra X = (\A^n)^r / S_n \ra A = \Spec R_\Sigma
\eeq
where the intervening quotient is just the affine scheme associated to
the ring of multisymmetric functions.  Every $r$-tuple of integers in
$\OO_L$ corresponds to an integral point of $X$; however, the fact
that $X$ fails to embed naturally in a low-dimensional affine space
makes it difficult to count points of $X(\Z)$ with bounded height.
The method used here identifies a locus $W \subset X$ which is
contracted in the map to $\Spec R_\Sigma$, and shows that the map $X(\Z) \ra
A(\Z)$ has fibers of bounded size away from $W$; this gives an upper
bound on the number of integral points on $X \backslash W$ of bounded
height.  One is naturally led to ask whether Manin's conjecture on
points of bounded height on varieties can be used to predict an
asymptotic number of integral points on some open subscheme of $X$.  Any
such prediction would lead to a refinement of our upper bound on the
number of number fields.

\subsection{Improvements, Invariant Theory, and the Large Sieve}
\begin{rem}
The method we have used above may be optimized
in various ways: by utilizing more of the invariant theory
of $S_n$, and by using results about counting integral points on varieties.
These techniques may be used, for any fixed $n$, to improve
the exponent in the upper bound of Theorem \ref{th:main}.
(The invariant theory, however, becomes more
computationally demanding as $n$ increases).
However, they do not change the limiting behavior as $n \rightarrow
\infty$.
We have therefore chosen to present a different example
of this optimization:
giving good bounds on $N_{K,n}(X;G)$ for $G \neq S_n$.
For simplicity of exposition we take $K = \Q$. 
\end{rem}

\begin{exmp}
Let $G = 
\langle (1,6,2)(3,4,5), (5,6)(3,4) \rangle$;
it is a primitive permutation group on $\{1,2,3,4,5,6\}$
whose action is conjugate to the action of $\mathrm{PSL}_2(\mathbb{F}_5)$
on $\mathbb{P}^1(\mathbb{F}_5)$. 

We will show $N_{\Q,6}(X;G) \ll_{\epsilon} X^{8/5 + \epsilon}$, a
considerable improvement over Schmidt's bound of $X^2$ (over which, in
turn, Theorem~\ref{th:main} presents no improvement for $n=6$).

Let $G$ act on monomials $x_1, x_2, \dots, x_6$ by permutation of the
indices.
Set $f_i = \sum_{j=1}^{6} x_j^i$ for $1 \leq i \leq 5$,
and $f_6 = x_1 x_2 (x_3 + x_4) + x_1 x_3 x_5 + x_1 x_4 x_6 + x_1 x_5 x_6 +
x_2 x_3 x_6 + x_2 x_4 x_5 + x_2 x_5 x_6 + x_3 x_4 (x_5 + x_6)$.
Set $A = \C[f_1,f_2,\dots,f_6]$.
Then $R = \C[x_1, \dots, x_6]^G$ is a free $A$-module of degree $6$;
indeed $R = \oplus_{i=1}^{6} A \cdot g_i$,
where $g_1 = 1$ and $g_2, g_3, \dots, g_5$ can be chosen to be homogeneous
of degree $5,6,6,7,12$. (This data was obtained
with the commands {\tt InvariantRing},{\tt PrimaryInvariants},
and {\tt SecondaryInvariants} in Magma.)
One checks that $\overline{R} = R/f_1 R$ is an integral domain. 
%Moreover $V = \C^6$, on which $G$ acts
%by coordinate permutations, decomposes into an invariant
%line $L$ ($x_1 = \dots = x_6$) and an invariant
%$5$-dimensional subspace $W$ ($\sum_{i=1}^{6} x_i = f_1 = 0$) under the
%action of 
%$G$; this gives a decomposition $V = L \oplus W$
%and correspondingly $\C[V] = \C[L] \otimes_{\C} \C[W]$,
%so $R = \C[L] \otimes_{\C} \C[W]^G$.
%From this one deduces that $\overline{R} = R/f_1 R$ is
%isomorphic to $\C[W]^G$, and therefore an integral domain.

Let $S$ be the subring of $\overline{R}$ generated
by  $\overline{f_2}, \dots, \overline{f_6}$ and $
\overline{g_2}$, and let 
$Z = \mathrm{Spec}(S)$.
$S$ is an integral domain since $\overline{R}$ is;
thus $Z$ is irreducible. The map
$\C[f_2,f_3,f_4,f_5,f_6] \rightarrow S$
induces a finite projection $Z \stackrel{\Pi}{\rightarrow} \mathbb{A}^5$
(it is finite since $R$ is finite over $A$, so $\overline{R}$
is finite over $\C[f_2,f_3, \dots, f_6]$). 
Also $\overline{g_2} \notin \C[\overline{f_2}, \dots,
\overline{f_6}]$, as follows from the fact that
$R = \oplus_{i=1}^{6} A g_i$;
thus the degree of $\Pi$ is at least $2$.

Suppose $L$ is a number field with $[L:\Q] = 6$
with Galois group $G$
and $\DD_{L} < X$.
Minkowski's theorem implies there exists $x \in \OO_L$
with $\Tr_{\Q}^{L}(x) = 0$ and $||x||\ll X^{1/10}$;
here $||x||$ is defined as in the proof of Proposition \ref{prop:rtuple}. 
The element $x \in \OO_L$
gives rise to a point $\mathbf{x} \in Z(\Z)$
whose projection $\Pi(\mathbf{x}) = (y_1,y_2,y_3,y_4,y_5) \in \Z^5$
satisfies:
\begin{equation} \label{eqn:skewbox}|y_1| \ll X^{2/10}, 
|y_2| \ll X^{3/10},
|y_3| \ll X^{4/10}, y_4 \ll X^{5/10}, y_5 \ll X^{3/10}
\end{equation}
We must count integral points on $Z$ whose projection
to $\mathbb{A}^5$ belong to
the skew-shaped box
defined by (\ref{eqn:skewbox}).
It is clear that the number of points on
$Z(\Z)$ projecting to the box (\ref{eqn:skewbox})
is at most $X^{17/10}$, but applying the large sieve
to the map $Z \stackrel{\Pi}{\rightarrow} \mathbb{A}^5$ (c.f.
\cite{sdcohen} or \cite{serre}) one obtains
the improved bound $X^{8/5+\epsilon}$. (Note that
the results, for example in \cite{serre},
are stated only for a ``square''  box
(all sides equal) around the origin --
but indeed they apply, with uniform implicit constant,
to a square box centered
at {\em any} point. Now we tile the
skew box (\ref{eqn:skewbox}) by
square boxes of side length $X^{2/10}$ to obtain the claimed
result.)

One expects that one can quite considerably improve this bound
given more explicit understanding of the variety $Z$; ideally
speaking one would like to slice it, show that most slices
are geometrically irreducible, and apply the Bombieri-Pila bound
\cite{pila}. It is the intermediate step -- showing
that very few slices have irreducible components of low degree --
which is difficult. This seems like an interesting computational question. 

We remark that this particular example can also be analyzed by constructing
an associated {\em quintic} extension (using the isomorphism of
$\mathrm{PSL}_2(\mathbb{F}_5)$ with $A_5$) and counting these quintic
extensions. This is close in spirit to the idea of the next section; in any
case the method outlined above should work more generally. 
\end{exmp}

\subsection{Counting Galois extensions}

In this section, we give bounds on the number of Galois extensions of
$\Q$ with bounded discriminant.  In combination with the lower bound
in Theorem~\ref{th:main} for the total number of extensions, this
yields the fact that ``most number fields, counted by discriminant,
are not Galois.''

Let $K$ be a number field of degree $d$ over $\Q$ and $G$ a finite
group; we denote by $N_K(X,G)$ the number of Galois extensions of $K$
with Galois group $G$ such that
\beq
\N^K_\Q \DD_{L/K} < X.
\eeq
 
\begin{prop} If $|G| > 4$, then $N_K(X,G) \ll_{K,G,\epsilon} X^{3/8+\epsilon}$.
\label{pr:galois}
\end{prop}

\begin{rem} Proposition~\ref{pr:galois} is not meant to be sharp; our
  aim here is merely to show that most fields are not Galois, so we
  satisfy ourselves with giving a bound smaller than $X^{1/2}$.  In
  fact, according to a conjecture of Malle~\cite{mall:jnt}
, $N_K(X,G)$ should be
  bounded between $X^{\frac{\ell}{(\ell - 1)|G|}}$ and $X^{
 \frac{\ell}{(\ell - 1)|G|}
  + \epsilon}$, where $\ell$ is the smallest prime divisor of $|G|$.
  This conjecture is true for all abelian groups $G$ by a theorem of
  Wright~\cite{wrig:abelian}, and is proved for all nilpotent groups
  in a preprint of Kl\"{u}ners and Malle~\cite{klun:nilpotent}.
\end{rem}

\begin{proof} We proceed by induction on $|G|$.  In this proof,
all implicit constants in $\ll, \gg$ depend on $K$,
$\epsilon$ and $G$, although we do not always explicitly note this. 

Write an exact sequence
\beq
1 \ra H \ra G \ra Q \ra 1
\eeq
where $H$ is a minimal normal subgroup of $G$.  Then $H$ is a direct
sum of copies of some simple group~\cite[3.3.15]{robi:groups}.

If $L/K$ is a Galois extension with
$G_{L/K} \cong G$ and $\N^K_\Q \DD_{L/K} < X$, and $M$ is the subfield
of $L$ fixed by $H$, then $M/K$ is a Galois extension with Galois
group $Q$ and $N^K_\Q \DD_{M/K} < X^{1/|H|}$.  The number of such
extensions $M/K$ is $N_K(X^{1/|H|}, Q)$, which by the induction hypothesis
is $\ll_{Q} X^{3/8|H| + \epsilon}$ in case $|Q| > 4$.  If $|Q| \leq 4$, then
$Q$ is abelian and by Wright's theorem $N_K(X^{1/|H|}, Q) \ll
X^{1/|H| + \epsilon}$.

Now take $M/K$ to be fixed; then the number of choices for $L$ is
bounded above by $N_M(X,H)$. 

First, suppose $H$ is not abelian.  Let $H_0$ 
be a proper subgroup of $H$ that does not contain
any normal subgroups of $H$ and is of maximal cardinality subject to this
restriction. 

If $L_1$ and $L_2$ are two Galois $H$-extensions of $M$, then $L_1 \cong
L_2$ if and only if $L_1^{H_0} \cong L_2^{H_0}$.  If $L' = L^{H_0}$,
then $\DD_{L'}  \leq (\DD_{L})^{1/|H_0|}
\ll_{K} (\N^{K}_{\Q} \DD_{L/K})^{1/|H_0|}$. 
It follows from the main theorem of \cite{schm:columbia}
that the number of possibilities for $L'$
(and hence the number of possibilities for $L$),
{\em given} $M$, 
is $\ll X^{\frac{(|H|/|H_0| + 2)}{4 |H_0|}}$,
where the implicit constant is independent of $M$. 
%we also have
%\beq
%\N^M_\Q \DD_{L/M} \geq (\N^M_\Q \DD_{L'/M})^{|H_0|}.
%\eeq
%It follows that
%\beq
%N_M(X,H) \leq N_{M,|H|/|H_0|}(X^{1/|H_0|}).
%\eeq
%Then by the main theorem of \cite{schm:columbia}, we have
%\beq
%N_M(X,H) \ll_{d|Q|,[H:H_0]} X^{(|H|/|H_0|+2)/4|H_0|}
%\eeq
%for some absolute constant $C$. 

The group $H_0$ can be chosen to have
size at least $\sqrt{|H|}$ (\cite{kl}, comments after 5.2.7)  so
%\beq
%N_M(X,H) \ll_{d|Q|,[H:H_0]} X^{1/4 + 1/(2\sqrt{|H|})}.
%\eeq
 summing over all choices of $M$, we find
\beq
N_K(X,G) \ll_{G} X^{1/4 + 1/2\sqrt{|H|} + 1/|H| + \epsilon}
\eeq
which, since $|H| \geq 60$, proves Proposition~\ref{pr:galois} in case $H$
is non-abelian.
 
Now, suppose $H$ is abelian; so $H = (\Z/p\Z)^r$ for some prime $p$
and some positive integer $r$.  By \cite{wrig:abelian} we may assume
$|Q| \geq 2$.

Let $b_M(Y)$ be the number of $H$-extensions of $M$ such that $\N^M_\Q
\DD_{L/M} = Y$.  Let $S$ be the set of primes of $\Q$ dividing $Y$,
let $G_S(M)$ be the Galois group of the maximal extension of $M$ unramified
away from primes dividing $S$, and for each prime $\lambda$ of $M$ let
$I_\lambda$ be the inertia group at $\lambda$.  Then $b_M(Y) \leq
|\Hom(G_S(M),H)|.$ Moreover, the kernel of the map
\beq
\Hom(G_S(M),H) \ra \bigoplus_{\lambda | S} \Hom(I_\lambda, H)
\eeq
is isomorphic to a subgroup of the $r$th power of the class group of
$M$, and as such has cardinality $\ll_\epsilon
\DD_{M/\Q}^{r/2+\epsilon}$, by the easy part of the Brauer-Siegel
theorem.  On the other hand, the number of primes $\lambda$ is $\ll
|S|$, and $|\Hom(I_\lambda, H)|$ is bounded by some constant $C$
depending only on $[M:\Q]$; so the
image of the map above has cardinality at most
\beq
(C')^{|S|} \ll_{\epsilon,K,G} Y^\epsilon.
\eeq
We conclude that
\begin{equation}
b_M(Y) \ll \DD_{M/\Q}^{r/2 + \epsilon} Y^\epsilon.
\label{eq:bmy}
\end{equation}

Let $\mu$ be a prime of $K$ such that $\mu$ does not divide $|G|
\DD_{M/K}$ and primes of
$M$ above $\mu$ ramify in $L$.  Then the image of $I_\mu \subset
\GalK$ in $G$ is a cyclic subgroup whose order is a multiple of $p$;
it follows that $(p-1)|G|/p$  divides $\ord_\mu \DD_{L/K}$.  So $\N^M_\Q
\DD_{L/M}$ lies in one of a finite set of cosets of
$\Q^*/(\Q^*)^{(p-1)|G|/p}$.
Let $\Sigma$ be this union of cosets.
Since the valuation of 
$\N^M_{\Q} \DD_{L/M}$ is divisible by $\frac{(p-1)|G|}{p}$ at primes not dividing $|G| \N^{K}_{\Q}
\DD_{M/K}$,
it follows that we may take $\Sigma$ so that
the number of cosets in $\Sigma$ is $\ll_{\epsilon,G} ( \N^{K}_{\Q} \DD_{M/K})^{\epsilon}$. 

When $M$ is a $Q$-extension of $K$, we write $N_1$ for $\N^K_\Q
\DD_{M/K}$.  Then
\beq
N_K(X,G) \leq \sum_{M:N_1 \leq X^{1/|H|}}
\sum_{{N_2 < X N_1^{-|H|},N_2 \in \Sigma}} b_M(N_2).
\eeq

%\begin{eqnarray*}
%N_K(X,G) & = & \sum_{N_1^{|H|} N_2 < X} a(N_1,N_2) \\
%& = & \sum_{N_1 < X^{1/|H|} \sum_{N_2 < X N_1^{-|H|} a(N_1,N_2) \\
%& \leq  & \sum_{N_1 < X^{1/|H|} \sum_{M:N^K_\Q \DD_{M/K} = N_1}
%\sum_{N_2 < X N_1^{-|H|}, N_2 \in \Sigma} b_M(N_2).
%\end{eqnarray*}
The inner sum has length $\ll_{\epsilon} N_1^{\epsilon}(X N_1^{-|H|})^{p/(p-1)|G|}$, which,
combined with \eqref{eq:bmy}, gives
\begin{eqnarray*}
N_K(X,G) & \ll_\epsilon & \sum_{M: N_1 \leq X^{1/|H|}}  X^{p/(p-1)|G|+ \epsilon} N_1^{r/2 - p/(p-1)|Q| + \epsilon} \\
& \leq & 
N_K(X^{1/|H|},Q) X^{p/(p-1)|G| + \epsilon} \max_{N_1 < X^{1/|H|}}
    N_1^{r/2 - p/(p-1)|Q| + \epsilon} \\
& = & N_K(X^{1/|H|},Q) X^{\alpha + \epsilon}
\end{eqnarray*}
where $\alpha = \max(\frac{r}{2|H|},\frac{p}{(p-1)|G|})$.

By the induction hypothesis, $N_K(X^{1/|H|},Q) \ll X^{3/8|H| + \epsilon}$ when $|Q|
\geq 5$, while $N_K(X^{1/|H|},Q)$ is asymptotic to
$X^{1/2|H|}$ if $|Q| = 3,4$ and to $X^{1/|H|}$ when $|Q| = 2$.  Define
\beq
\beta(Q) = \left\{\begin{array}{ll}
3/8 & |Q| \geq 5; \\
1/2 & |Q| =3,4; \\
1 & |Q| = 2 \\
\end{array}\right.
\eeq
Then 
\beq
N_K(X^{1/|H|},Q) X^{r/2|H|} \ll
X^{(r/2 + \beta)/|H| + \epsilon}
\eeq
and the exponent $\frac{r/2 + \beta}{|H|}$ is at most $3/8$ unless either
$|H| = 2$, or $|Q| = 2$ and $|H| = 3,4$.  In case $|Q| = 2, |H| = 4$,
the group $G$ is nilpotent and Proposition~\ref{pr:galois} is proved
by Kl\"{u}ners and Malle. 

%If $|H| = 4$, the
%case $|Q| = 2$ can be neglected, since in that case $G$ is nilpotent
%and thus covered by the result of Kl\"{u}ners and Malle.  So $(r/2 +
%\beta)/|H| \leq (1+1/2)/4 = 3/8$.  If $|H| = 3$, then $(r/2 +
%\beta)/|H| < 3/8$ unless $|Q| = 2$.
On the other hand,
\beq
N_K(X^{1/|H|},Q) X^{p/(p-1)|G|} \ll
X^{(p/(p-1)|Q| + \beta)/|H| + \epsilon}.
\eeq
Here, the exponent is once again at most $3/8$ unless either $|H| =
2$, or $|Q| = 2$ and $|H| = 3,4$.

We have thus proven Proposition~\ref{pr:galois} unless $G = S_3$ or $H
= \Z/2\Z$.  In the former case, the proposition follows from the
theorem of Datskovsky and Wright~\cite{dats:cubics} on the number of
cubic extensions of number fields.  In the latter case, we can refine
the argument above; let $b'_M(Y)$ be the number of quadratic
extensions $L/M$ which are preserved by the action of $Q$
and so that $\N^{M}_{\Q} \DD_{L/M} = Y$. Choosing 
$S$ to consist of all divisors of $Y \N^{K}_{\Q} \DD_{M/K}$
and utilising the
inflation-restriction sequence
\beq
\Hom(G_S(K),\Z/2\Z) \ra \Hom(G_S(M),\Z/2\Z)^Q \ra H^2(Q,\Z/2\Z)
\eeq
we see that $b'_M(Y) \ll (Y \cdot \N^{K}_{\Q}\DD_{M/K})^{\epsilon}$. (Here $G_S(K)$ is defined
analogously to $G_S(M)$.) This saves a factor of
  $N_1^{r/2}$ throughout the rest of the argument, and in particular
  we have
\beq
N_K(X,G) \ll_\epsilon N_K(X^{1/2},Q) X^{2/|G| + \epsilon}.
\eeq
Since we may assume $G$ non-nilpotent, we can take $|Q| \geq 6$, which
yields
\beq
N_K(X,G) \ll_\epsilon X^{3/16} X^{1/6 + \epsilon}
\eeq
which again yields the desired result.
\end{proof}

\section{Proof of Lower Bound for $S_n$ extensions}

We now turn to the (easier) question of proving the
lower bounds for $N_{K,n}'(X)$ asserted in Theorem \ref{th:main},
 and finish with a brief discussion of some related issues. 

We make some preliminary remarks. 
Firstly, as was discussed in the Introduction,
this question is often much easier if one is counting extensions for 
$G$ a {\em proper} subgroup of $S_n$
(see Malle \cite{mall:jnt} for some examples). 
On the other hand, the general
question of lower bounds subsumes the inverse Galois problem over $\Q$.
The method we give can be generalized to $G$-extensions,
so long as one can construct a family of polynomials with generic Galois
group $G$.

As before, let $K$ be a fixed extension of $\Q$ of degree $d$.
%If $L$ is an extension of $K$ of degree $n$,
%and $\gamma \in \OO_L$, we define $|\gamma|_2$ 
%as in the proof of Proposition \ref{prop:rtuple}. 
We also set $\Delta_L = \N^{K}_{\Q}(\DD_{L/K})$
and $\OO_L^{0} = \{x \in \OO_L: \Tr^L_K(x) = 0\}$. 
In this section, we will not aim for any uniformity in $K$;
the implicit constants in this section will always depend on $K$ and $n$. 
As before, for $x$ an algebraic integer, we denote by $||x||$ the largest
archimedean valuation of $x$.

\begin{lem} \label{lemma:ceven}
Let $[L:K] = n$ be so that $L/K$ has no proper subextensions. 
Then $||x|| \gg \Delta_L^{\frac{1}{n(n-1)d}}$ for all $x \in \OO_L^{0}$.
\end{lem}

\proof 
If $x \in \OO_L^{0}$, then $\OO_K[x]$ is a subring of $\OO_L$;
in particular, the discriminant $\mathcal{D}(\OO_K[x])$ 
of $\OO_K[x]$ over $\OO_K$
is divisible
by $\DD_{L/K}$. In particular, $\N^K_{\Q}(\DD_{L/K})
\leq \N^K_{\Q} \mathcal{D}(\OO_K[x])$. 
$\mathcal{D}(\OO_K[x])$ 
is the same as the discriminant of the minimal polynomial
of $x$; from this, one deduces that $\DD(\OO_K[x])$ is
a principal ideal of $\OO_K$, generated by a polynomial of degree
$n(n-1)$ in the Galois conjugates of $x$.   
In particular, one deduces $\N^{K}_{\Q}(\mathcal{D}(\OO_K[x]))
\ll ||x||^{n(n-1) d}$, whence the assertion. 
\qed

See Remark \ref{lattices} for generalizations.

In the lower bound proved below,
we have not aimed to optimize the exponent
$1/2 + 1/n^2$. 
It will be obvious
from the proof that it can be improved somewhat, both
by replacing Schmidt's upper bound with that of Theorem
\ref{th:main}, and by utilizing successive maxima
and Remark \ref{lattices}
rather than just Lemma \ref{lemma:ceven}.
This seems like an interesting optimization question;
the gain for small $n$ can be significant although one does not obtain
an exponent near $1$. 
\proof (of lower bound 
$N_{K,n}'(X) \gg_{K,n} X^{1/2 + 1/n^2}$ in Theorem \ref{th:main}). 

We fix as before an algebraic closure $\bar{K}$.
Consider the set $S(Y)$ of algebraic integers $x \in \bar{K}$ so that
$[K(x):K] = n$,
$\Tr^{K(x)}_{K}(x) = 0$
and $||x|| \leq Y$.  Let $S(Y;S_n)$ be the subset
of those $x$ so that the Galois closure of $K(x)$ over $K$
has Galois group $S_n$.

Then, by considering the characteristic polynomial, we see that
$|S(Y)| \gg Y^{d(n(n+1)/2-1)}$. Considering (the proof of)
Hilbert's irreducibility theorem, we see that the same bound holds 
for $S(Y;S_n) \subset S(Y)$:
\begin{equation}
\label{eqn:hi}|S(Y;S_n)| \gg Y^{d(\frac{n(n+1)}{2} -1)}=Y^{\frac{(n-1)(n+2)d}{2}}
\end{equation}
Indeed one may put a congruence constraint on the characteristic
polynomial to guarantee that the Galois closure has group $S_n$ (c.f. \cite{schinzel}).

Suppose $L$ is an $S_n$-extension of $K$ (i.e. $[L:K]  = n$
and the Galois closure of $L/K$ has Galois group $S_n$).
$\OO_L^{0}$ is a free $\Z$-module of rank $(n-1)d$; then
Lemma \ref{lemma:ceven} guarantees
that the number of $x \in S(Y;S_n)$ such that $K(x) \cong L$
is $ \ll (\frac{Y}{\Delta_L^{1/n(n-1)d}})^{(n-1)d}$;
in particular if there is at least one such $x$,
one must have $\Delta_L \ll Y^{n(n-1)d}$. 
Combining these comments with (\ref{eqn:hi}) we find
that for some constant $c$:
\begin{equation} \label{temp}\sum_{
\stackrel{L: \Delta_L \leq c Y^{n(n-1)d}}{L/K \, S_n-extension}}
\left(\frac{1}{\Delta_L}\right)^{\frac{1}{n}} \gg
Y^{dn(n-1)/2}\end{equation}
However, Schmidt's upper bound $N_{K,n}(X) \ll X^{(n+2)/4}$
easily shows that $$\sum_{\stackrel{L: \Delta_L < Y^{d(n-1)}}{[L:K]=n}}
\left(\frac{1}{\Delta_L}\right)^{\frac{1}{n}} \ll Y^{\frac{d n(n-1)}{2} -
\delta}$$ for some
$\delta > 0$; thus one can replace the range of summation
in (\ref{temp}) by $Y^{d(n-1)} < \Delta_L \leq c Y^{d n(n-1)}$ without changing
the result. In particular $N_{K,n}'(c Y^{d n(n-1)})
\gg Y^{d n(n-1)(\frac{1}{2} + \frac{1}{n^2})} $, which implies the result. 
\qed

\begin{rem} \label{lattices} (Shape of number field lattices)
Lemma \ref{lemma:ceven} emphasizes the importance of understanding
the shape of number field lattices. 
For clarity, 
fix attention on totally real number fields of degree $n \geq 3$
over $\Q$
with no proper subfields;
one can formulate similar ideas in the general case.

Let $L$ be such a number field. Then $\OO_L^{0}$ is 
a lattice endowed with a natural quadratic form, 
namely $x 
\mapsto \mathrm{tr}(x^2)$; as such,
it defines an element $[L]$ of the moduli space $\Siegel$
of {\em homothety classes of positive definite quadratic forms}.
$\Siegel$
can be identified with $\PGL_{n-1}(\Z) \backslash \PGL_{n-1}(\R)/PO_{n-1}(\R)$. 
It is reasonable to ask about the distribution of $[L]$,
as $L$ varies,
in the finite volume space $\Siegel$. 

Hendrik Lenstra has informed us that David Terr has proven the equidistribution of a closely related
set in the case $n=3$ in his Ph.D. thesis.

General results in this direction seem out of reach; one can at least prove, however, mild constraints on $[L]$
that show it does not lie {\em too} far into the cusp. 
Let $a_1\leq a_2 \leq \dots \leq a_{n-1}$ be the successive minima 
(in the sense of Minkowski) of 
$\OO_L^{0}$. Then one has automatically
$a_1 a_2 \dots a_{n-1} \asymp \sqrt{\Disc_{L}}$;
however, on account of the assumption that $K$ has no proper subfield, one further has
for $1 \leq j \leq n-2$ that
$a_1 a_j \gg a_{j+1}$ (indeed, were this not so, 
the lattice spanned by $a_1, a_2, \dots, a_j$
would be stable under multiplication by $a_1$,
and so $\Q(a_1)$ is a proper subfield of $L$). Finally
one evidently has $a_j \gg 1$. Combining
these constraints gives nontrivial constraints on the $a_i$;
for example, one recovers Lemma \ref{lemma:ceven},
and one obtains $a_{n-1} \ll  \Disc_{L}^{\frac{1}{2 ([\sqrt{2n}] -
1)}}$, where $[\alpha]$ is the greatest integer $\leq \alpha$. 
One may use this type of result to further improve the exponents in Theorem 
\ref{th:main}
for specific $n$. 
\end{rem}

\begin{rem} (Alternate ways of ordering number fields)
There are many ways to order lattices of rank $>1$;
the ordering by volume is completely different than that by shortest
vector. 

We continue to work over the base field $\Q$. 
Given a number field $L$, we define $\short(L) = 
\mathrm{inf}(||x||: x \in \OO_L, \Q(x) = L)$.
It is then immediate that, for any $C > 0$, the number of number fields $L$
with $[L:\Q] = n$ and $\short(L) \leq C$ is finite; indeed one may verify
that $\short(L)$ is ``comparable'' to the discriminant:
$\Disc_{L}^{\frac{1}{n(n-1)}} \ll \short(K)  \ll
\Disc_{L}^{\frac{1}{2 [(n-1)/2]}}$.

Let $N_{n,\short}(Y)$ be the number of $L$ with $[L:\Q] = n$ and 
$\short(L) \leq Y$. Then one may show quite easily that
$Y^{\frac{(n-1)n}{2}} \ll N_{n,\short}(Y) \ll Y^{\frac{(n-1)(n+2)}{2}}$;
in particular, the discrepancy between upper and lower bounds is much
better than when counting by discriminant. 
Further, the (approximate)
asymptotic $N_{n,\short}(Y) \asymp Y^{\frac{(n-1)(n+2)}{2}}$
follows from the Hypothesis below, which seems
very difficult (Granville \cite{granville}
and Poonen have proved versions of this -- too weak for our purposes --
using the $ABC$ conjecture.) The idea
is to use Hypothesis \ref{star} to construct many polynomials with square-free
discriminant. 

\begin{hyp}\label{star}
Let $f \in \Z[x_1, \dots, x_n]$. Then, if $B_i$ is any sequence of boxes
all of whose side lengths go to infinity, one has:
$$\lim_{i} \frac{\#\{x \in B_i: f(x) \mbox{ squarefree}\}}{\#\{x \in B_i\}} =
C_f$$
where $C_f$ is an appropriate product of local densities.
\end{hyp}

\end{rem}

\bibliographystyle{plain}
\bibliography{CountingNFf}

%\begin{thebibliography}{10}

%\bibitem{bhar:thesis}
%M. Bhargava
%\newblock {\em Higher Composition Laws}.
%\newblock Ph.D. thesis, Princeton University, 2001.
 
%\bibitem{bhar:hcl4}
%M. Bhargava
%\newblock Higher composition laws, IV.
%\newblock Preprint.

%\bibitem{dats:cubics}
%B. Datskovsky and D. J. Wright.
%\newblock  Density of discriminants of cubic extensions.
%\newblock {\em  J. Reine Angew. Math.}  386  (1988), 116--138.
%
%\bibitem{kabl:quintics}
%A. Kable and A. Yukie
%\newblock On the number of quintic fields.
%j\newblock Preprint.
%
%\bibitem{katz:mgpoints}
%A.J. de Jong and N.M.Katz.
%\newblock Personal communication.
%
%\bibitem{klun:nilpotent}
%J. Kl\"{u}ners and G. Malle.
%\newblock Counting nilpotent Galois extensions.
%\newblock Preprint.
%
%\bibitem{robi:groups}
%D.J.S. Robinson.
%\newblock {\em A Course in the Theory of Groups.}
%\newblock Springer-Verlag, Berlin, 1993.
%
%\bibitem{sieg:gon}
%C.L.Siegel.
%\newblock {\em Lectures on the Geometry of Numbers.}
%\newblock Springer-Verlag, Berlin, 1989.
%
%\bibitem{yuki:quartics}
%A. Yukie.
%\newblock Density theorems related to prehomogeneous vector spaces.
%\newblock Preprint.
%
%\end{thebibliography}

\end{document}